\documentclass[11pt]{article}

\usepackage{amsmath}
\usepackage{amssymb}
\usepackage{eucal}

\begin{document}

\baselineskip 16pt

\title{A Robinson characterization of finite $P\sigma T$-groups }

\author            
{ Alexander  N. Skiba \\
{\small  Department of Mathematics and Technologies of Programming, 
  Francisk Skorina Gomel State University,}\\
{\small Gomel 246019, Belarus}\\
{\small E-mail: alexander.skiba49@gmail.com}}

\date{}
\maketitle

\begin{abstract}     Let   $\sigma =\{\sigma_{i} | i\in I\}$ be some
 partition of the set of all primes $\Bbb{P}$ and  let $G$ be a  finite group. 
Then  $G$
 is said to be   \emph{$\sigma $-full} if $G$ has a Hall $\sigma _{i}$-subgroup for all $i$.  
A subgroup $A$ of $G$ is said to be    \emph{$\sigma$-permutable} in $G$ provided 
 $G$ is  $\sigma $-full and 
$A$ permutes with all Hall $\sigma _{i}$-subgroups $H$ of $G$ (that is, $AH=HA$)  for all $i$.  

   We obtain a  characterization of    finite  groups $G$ 
 in which  $\sigma$-permutability    is a transitive relation in $G$,
 that is, if $K$ is a ${\sigma}$-permutable subgroup
 of $H$ and 
 $H$ is a ${\sigma}$-permutable subgroup of $G$, then  $K$ is a
 ${\sigma}$-permutable subgroup of $G$.

\end{abstract}


\let\thefootnote\thefootnoteorig

\section{Introduction}

Throughout this paper, all groups are finite and $G$ always denotes
a finite group. Moreover,  $\mathbb{P}$ is the set of all  primes,
  $\pi= \{p_{1},  p_{2}, \ldots \} \subseteq  \Bbb{P}$ and  $\pi' =  \Bbb{P} \setminus \pi$.
 If
 $n$ is an integer, the symbol $\pi (n)$ denotes
 the set of all primes dividing $n$; as usual,  $\pi (G)=\pi (|G|)$, the set of all
  primes dividing the order of $G$.

If  $1\in \mathfrak{F}$ is a class of groups, then  $G^{\mathfrak{F}}$ 
 denotes the \emph{$ \mathfrak{F}$-residual of $G$}, that is, 
intersection of all normal subgroups $N$ of $G$ with  $G/N\in \mathfrak{F}$;
 $G_{\mathfrak{F}}$ 
 denotes the \emph{$ \mathfrak{F}$-radical of $G$}, that is, 
the product  of all normal subgroups $N$ of $G$ with  $N\in \mathfrak{F}$.

In what follows, $\sigma$  is some partition of  
$\Bbb{P}$, that is,  $\sigma =\{\sigma_{i} |
 i\in I \}$, where   $\Bbb{P}=\bigcup_{i\in I} \sigma_{i}$
 and $\sigma_{i}\cap
\sigma_{j}= \emptyset  $ for all $i\ne j$; $G$
 is said to be   \emph{$\sigma $-full} \cite{2, commun} if $G$ has a
 Hall $\sigma _{i}$-subgroup for all $i$.

{\bf Definition 1.1.} We say that a subgroup $A$ of $G$ is   
   \emph{$\sigma$-permutable} in $G$ \cite{1} provided  $G$ is  $\sigma $-full and 
$H$ permutes with all Hall $\sigma _{i}$-subgroups $H$ of $G$ (that is, $AH=HA$)  for all $i$.

{\bf Remark 1.2.}   A     set  ${\cal H}$ of subgroups of $G$ is a
 \emph{complete Hall $\sigma $-set} of $G$ \cite{2, commun}  if
 every member $\ne 1$ of  ${\cal H}$ is a Hall $\sigma _{i}$-subgroup of $G$
 for some $\sigma _{i} \in \sigma$ and ${\cal H}$ contains exactly one Hall
 $\sigma _{i}$-subgroup of $G$ for every  $i$. 
By   Proposition 3.1 in \cite{???}, a subgroup $A$ of  $G$ is   
   $\sigma$-permutable  in $G$ if and only if  $G$ 
 possesses a 
 a complete Hall $\sigma$-set  ${\cal H}$ such that $AL^{x}=L^{x}A$ for all
  $L\in {\cal H}$ and  all $x\in G$.

Recall that   $G$ is said to be:  
\emph{$\sigma$-primary}   \cite{1} if  $G$ is a $\sigma_{i}$-group for some 
$i$,  \emph{$\sigma$-decomposable} (Shemetkov \cite{Shem})  or 
\emph{$\sigma$-nilpotent}  (Guo and Skiba \cite{33}) if $G=G_{1}\times \dots \times G_{n}$ 
for some $\sigma$-primary groups $G_{1}, \ldots, G_{n}$.

The  usefulness   of $\sigma$-permutable subgroups is connected  mostly 
with the following   their  property.

{\bf Theorem  A.} (See Theorem B in \cite{1}).  {\sl If $A$
 is 
a $\sigma$-permutable subgroup of $G$, then   
 $A^{G}/A_{G}$ is $\sigma$-nilpotent.}

{\bf Example  1.3.} (i) In the classical case,
 when $\sigma =\sigma ^{0}=\{\{2\}, \{3\}, \ldots 
\}$, the subgroup $A$ of  $G$  is ${\sigma} ^{0}$-permutable in $G$ if and only
 if $A$ permutes with all Sylow subgroups of $G$.
 Note that a ${\sigma} ^{0}$-permutable subgroup
 is  also called    
 \emph{$S$-permutable} \cite{prod}.  Note also that for every  $S$-permutable 
subgroup $A$ of $G$ the quotient  $A^{G}/A_{G}$ is nilpotent (Kegel, Deskins)  by Theorem A.

(ii) In the other classical case, when  $\sigma =\sigma ^{\pi}=\{\pi, 
\pi'\}$, a subgroup $A$ of  $G$  is ${\sigma} ^{\pi}$-permutable in $G$ if and only
 if  $G$ has  a Hall $\pi$-subgroup and a Hall $\pi'$-subgroup  and 
  $A$ permutes  with 
all Hall $\pi$-subgroups  and with all Hall $\pi'$-subgroups of $G$. For every 
 $\sigma ^{\pi}$-permutable 
subgroup $A$ of $G$  the quotient  $A^{G}/A_{G}$ is $\pi$-decomposable, that is,
  $A^{G}/A_{G}=O_{\pi}(A^{G}/A_{G})\times O_{\pi'}(A^{G}/A_{G})$  by Theorem A.

(iii) In fact, in the theory of $\pi$-soluble groups ($\pi= \{p_{1}, \ldots , p_{n}\}$)
 we deal with the  partition 
$\sigma =\sigma ^{0\pi }=\{\{p_{1}\}, \ldots , \{p_{n}\}, \pi'\}$ of $\Bbb{P}$. 
 The subgroup $A$ of 
 $G$ is $\sigma^{0\pi }$-permutable in $G$ if and only if  $G$ has
 a Hall $\pi'$-subgroup   and 
 $A$ permutes  with 
all Hall $\pi'$-subgroups  and with all Sylow $p$-subgroups  of $G$ for all $p\in \pi$.
 For every 
 $\sigma ^{0\pi}$-permutable 
subgroup $A$ of $G$ the quotient   $A^{G}/A_{G}$ is $\pi$-nilpotent, that is,
 $A^{G}/A_{G}=O_{\pi}(F(A^{G}/A_{G}))\times O_{\pi'}(A^{G}/A_{G})$  by Theorem A.

We say, following \cite{1},  that  $G$ is a {\sl $P\sigma T$-group} 
 if ${\sigma}$-permutability  
is a transitive relation in $G$, that is, if $K$ is a ${\sigma}$-permutable subgroup
 of $H$ and 
 $H$ is a ${\sigma}$-permutable subgroup of $G$, then  $K$ is a
 ${\sigma}$-permutable subgroup of $G$.   
 In the case when $\sigma =\sigma ^{0}$, a  $P\sigma 
T$-group is  also called a \emph{$PST$-group} \cite{prod}. 

 Note that if 
$G=(Q_{8}\rtimes C_{3})\curlywedge(C_{7}\rtimes C_{3})$   (see \cite[p. 50]{hupp}), where 
$Q_{8}\rtimes C_{3}=SL(2, 3)$   and  $C_{7}\rtimes C_{3}$ is a non-abelian 
group of order 21, then $G$ is not a $PST$-group but $G$ is  a  $P\sigma T$-group,
 where $\sigma =\{\{2, 3\},  \{2, 3\}'\}$

The   description of
 $PST$-groups   
 was first obtained by   Agrawal \cite{Agr},
  for the soluble  case, and 
  by Robinson in \cite{217}, for the general case. 
In the   further publications,   authors (see, for example,   the recent papers
 \cite{78}--\cite{rend} and 
 Chapter 2 in \cite{prod}) have found out and  described 
many other   interesting characterizations   of  $PST$-groups.

  In the case when $G$
 is \emph{$\sigma$-soluble} (that is, every chief factor of $G$ is $\sigma$-primary) 
 the description  of $P\sigma T$-groups  was obtained in the paper
\cite{ariv} on the  base of the results  and methods in \cite{1, basel, B, W}.

{\bf Theorem B } (See  Theorem A in
 \cite{ariv}).   {\sl        If  $G$ is a  $\sigma$-soluble 
 $P\sigma 
T$-group and  $D=G^{\frak{N_{\sigma}}}$ is the $\sigma$-nilpotent residual of $G$,
 then  
    the following conditions hold:} 

(i) {\sl $G=D\rtimes M$, where $D$   is an abelian  Hall
 subgroup of $G$ of odd order, $M$ is $\sigma$-nilpotent  and  every element of $G$ induces a
 power automorphism in $D$;  }

(ii) {\sl  $O_{\sigma _{i}}(D)$ has 
a normal complement in a Hall $\sigma _{i}$-subgroup of $G$ for all $i$.}

{\sl Conversely, if  Conditions (i) and (ii) hold for  some subgroups $D$ and $M$ of
 $G$, then $G$ is  a $P\sigma 
T$-group.}
 
Before continuing, we give some further definitions.
                              
{\bf Definition 1.4.}  We say   that $G$ is:

(i) \emph{$\sigma $-supersoluble}   
if every chief factor of $G$ below $G^{{\mathfrak{N}}_{\sigma}}$ is 
cyclic;

(ii) a  \emph{$\sigma $-$SC$-group}
if every chief factor of $G$ below $G^{{\mathfrak{N}}_{\sigma}}$ is  simple.

{\bf Example 1.5.}  (i)   $G$ is supersoluble if and only if $G$ is $\sigma 
$-supersoluble  where $\sigma =\sigma ^{0}$ (see Example 1.3(i)). 

(ii) The  group $G$ is called an \emph{$SC$-group} (Robinson \cite{217}) or 
a \emph{$c$-supersoluble group} (Vedernikov \cite{veder}) if every chief 
factor of $G$ is a simple group.  Note that   $G$ is  an 
$SC$-group   if and only if $G$  is $\sigma$-$SC$-group 
 where $\sigma =\sigma ^{0}$.

(iii)     Let $G=A_{5}\times B$, where $A_{5}$ is
 the alternating group of degree 5 and $B=C_{29}\rtimes C_{7}$ is 
 a non-abelian group of order 203, and let $\sigma =\{\{7\},
 \{29\}, \{2, 3, 5\},  \{2, 3, 5,  7, 29\}'\}$. Then  $G^{{\frak{N}}_{\sigma}}=C_{29}$, so
   $G$ is a $\sigma$-supersoluble group but it  is neither soluble nor  $\sigma$-nilpotent.

(iv)      Let $G=SL(2, 7)\times A_{7}\times A_{5}\times B$, where  $B=C_{43}\rtimes C_{7}$
is  a non-abelian group of order 301,  and let 
 $\sigma =\{\{2, 3, 5\}, \{7, 43\}, \{2, 3, 5, , 7, 43\}'\}$.
 Then  $G^{{\frak{N}}_{\sigma}}=SL(2, 7)\times A_{7}$, so
   $G$ is a $\sigma $-$SC$-group  but it  is not a $\sigma$-supersoluble group.

In what follows,     ${\mathfrak{U}}_{\sigma}$ is the class of
 all $\sigma $-supersoluble groups; ${\mathfrak{U}}_{c\sigma}$  is 
 the class of all $\sigma $-$SC$-groups.
      
We say that $G$ is \emph{$\sigma$-perfect}  if $G^{{\mathfrak{N}}_{\sigma}}=G$, that is, 
$O^{\sigma _{i}}(G)=G$ for all $i$.

From Theorem B it follows that every  $\sigma $-soluble  $P\sigma T$-group is 
 $\sigma $-supersoluble. Our first observation shows that in general case
 every    $P\sigma T$-group is a   $\sigma$-$SC$-group.

{\bf Proposition  A.}   {\sl  Let $G$ be a $P\sigma T$-group  and let 
  $D=G^{{\mathfrak{S}}_{\sigma}}$ be the $\sigma$-soluble residual of $G$.
 Suppose 
that  $G$ possesses a complete Hall $\sigma$-set $\cal H$ whose members 
are $PST$-groups.  Then the following conditions 
hold:   }

(i) {\sl  $G$ is a   $\sigma$-$SC$-group. }

(ii) {\sl $D=G^{{\mathfrak{U}}_{\sigma}}$ is $\sigma$-perfect and $G/D$ is a $\sigma$-soluble $P\sigma T$-group.}

(iii) {\sl   $G$  satisfies 
 $N_{\sigma _{i}}$ for all $i$.}

In this proposition  we say that  \emph{$G$  satisfies 
 $N_{\sigma _{i}}$}  if whenever $N$ is  a $\sigma$-soluble normal 
subgroup of $G$, $\sigma _{i}'$-elements of $G$ induce power automorphisms in 
  $O_{\sigma _{i}}(G/N)$.  We say also,  following \cite[2.1.18]{prod}, that  \emph{$G$  satisfies 
 $N_{p}$}  if whenever $N$ is  a soluble normal 
subgroup of $G$, $p'$-elements of $G$ induce power automorphisms in 
  $O_{p}(G/N)$.

{\bf Corollary  1.6}   (See Proposition 2.1.1 in \cite{prod}). 
 {\sl  Let $G$ be a $PST$-group. Then: }

(i) {\sl  $G$ is an   $SC$-group, and  }

(ii) {\sl   $G$ of satisfies 
 $N_{p}$ for every prime $p$.}

 {\bf Definition 1.7.}   We say    that 
$(D, Z(D); U_{1},  \ldots , U_{k})$  is a \emph{Robinson $\sigma$-complex }
 (a \emph{Robinson complex}  in the case $\sigma =\sigma ^{0}$) of $G$
  if the following fold:

(i) $D$ is a  $\sigma$-perfect normal  subgroup of $G$,

(ii) $D/Z(D)=U_{1}/Z(D)\times \cdots \times U_{k}/Z(D)$, where $U_{i}/Z(D)$ is a
 non-abelian simple chief factor of $G$ for all $i$,

(iii) every chief factor of $G$ 
below $Z(D)$ is cyclic, and 

(iv) $D^{0}\leq D$ for every normal subgroup $D^{0}$ of $G$ satisfying  
Conditions (i), (ii) and (iii).

{\bf Example 1.8.}  Let $G=SL(2, 7)\times A_{7}\times A_{5}\times B$ be the group in Example 
1.5(iv) and 
 $\sigma =\{\{2, 3, 5\}, \{7, 43\}, \{2, 3, 5,  7, 43\}'\}$.
 Then $$(SL(2, 7)\times A_{7}, Z(SL(2, 7)); SL(2, 7), A_{7}Z(SL(2, 7)))$$
 is a Robinson $\sigma $-complex of $G$ and
 $$(SL(2, 7)\times A_{7}\times A_{5}, Z(SL(2, 7)); SL(2, 7), A_{7}Z(SL(2, 7)),
 A_{5}Z(SL(2, 7))))$$  is a  Robinson complex of $G$.

 Being based on Theorems A and B  and using some  
ideas in   \cite{217, basel, B, W}, in the given paper we prove the following

{\bf Theorem  C.}  {\sl Suppose 
that  $G$ possesses a complete Hall $\sigma$-set $\cal H$  whose members 
are $PST$-groups.
 Then $G$ is a $P\sigma T$-group if
 and 
only if  $G$ has a  $\sigma$-perfect normal  subgroup $D$ such that:
 }

(i) {\sl  $G/D$ is a $\sigma$-soluble $P\sigma T$-group. }

 (ii) {\sl If $D\ne 1$, then $G$ has a Robinson
 $\sigma$-complex of the form $(D,  Z(D); U_{1},  \ldots U_{k})$,
 and }

(iii)   {\sl   If  $\{i_{1}, \ldots , i_{r}\}\subseteq \{1, \ldots , k\}$, where
 $1\leq r  < k$, then $G$ and $G /U_{i_{1}}'\cdots U_{i_{r}}'$ satisfy
 $N_{\sigma _{i}}$ for all  $i$ such that $\sigma_{i}\cap \pi (Z(D))\ne 
 \emptyset  \}$. }

{\bf Corollary   1.9} (Robinson  \cite{217}).  {\sl A group    $G$ is a $PST$-group if
 and 
only if  $G$ has a perfect normal  subgroup $D$ such that:}

(i) {\sl  $G/D$ is a soluble $PST$-group. }

 (ii) {\sl If $D\ne 1$, then $G$ has a Robinson complex of the form
 $(D,  Z(D); U_{1},  \ldots U_{k})$,
 and }

(iii)   {\sl   If  $\{i_{1}, \ldots , i_{r}\}\subseteq \{1, \ldots , k\}$, where
 $1\leq r  < k$, then $G$ and $G /U_{i_{1}}'\cdots U_{i_{r}}'$ satisfy
 $N_{p}$ for all $p\in \pi(Z(D))$. }

The class $1\in \mathfrak{F}$ is said to be a \emph{formation} if every
 homomorphic image of $G/G^{\mathfrak{F}}$  belongs to $ \mathfrak{F}$ for every
 group $G$, that is, if $G\in \mathfrak{F}$, then also  every homomorphic image of $G$ belongs
 to $\mathfrak{F}$ and $G/N\cap R\in \mathfrak{F}$ whenever $G/N\in \mathfrak{F}$
 and $G/ R\in \mathfrak{F}$.
  The formation
$\mathfrak{F}$ is said to \emph{(normally) hereditary}
 if $H\in \mathfrak{F}$ whenever $ G \in {\cal
F}$   and $H$ is a (normal) subgroup of $G$.

 We  prove   Proposition  A and Theorem C in Section 3. But before, in 
Section 2,  we study properties of $\sigma$-supersoluble groups and 
$\sigma$-$SC$-groups. In particular, we prove the following two results.

{\bf Proposition  B.}   {\sl  For any partition $\sigma$ of $\Bbb{P}$ the 
following hold: }

(i) {\sl The class ${\mathfrak{U}}_{c\sigma}$ is a normally hereditary
 formation. }

(ii) {\sl The class ${\mathfrak{U}}_{\sigma}$ is a  hereditary
 formation. }

{\bf Theorem  D}   {\sl Let $N=
G^{{\mathfrak{N}}_{\sigma}}$ and let    $D=N^{\mathfrak{S}}$ be the soluble residual of $N$.
 Then $G$ is a $\sigma$-$SC$-group if  and 
only if the following hold:    }

(i)  {\sl $D=G^{\frak{U}_{\sigma }}$, and }

 (ii) {\sl if $D\ne 1$, then $G$ has   a 
Robinson complex  of the form $(D,  Z(D); U_{1}, \ldots , U_{k})$,
 where $Z(D)=D_{\mathfrak{S}}$ is  the soluble radical of $D$.}

{\bf Corollary  1.10}  (Robinson  \cite{217}).  {\sl  A group $G$ is an
 $SC$-group if  and 
only if $G$ satisfies:   }

(i) {\sl  $G/G^{\mathfrak{S}}$ is supersoluble.}

(ii) {\sl If $D=G^{\mathfrak{S}}\ne 1$, then  $G$ has   a 
Robinson complex of the form  $(D,  Z(D); U_{1}, \ldots , U_{k})$.}

 \section{Proofs of Proposition B and Theorem B}

The following lemma collects the properties of  $\sigma$-nilpotent groups
which we use in our proofs.

{\bf Lemma 2.1 } (See  Corollary 2.4 and Lemma 2.5  in \cite{1}).  {\sl  
The class   of all  $\sigma$-nilpotent groups
 ${\mathfrak{N}}_{\sigma}$  is closed under taking  
products of normal subgroups, homomorphic images and  subgroups. 
  }

{\bf Lemma 2.2} (See  \cite[2.2.8]{15}). {\sl If $\frak{F}$ is a  
formation and 
$N$, $R$ are   subgroups of $G$, where $N$ is normal in $G$, then}

(i) {\sl 
 $(G/N)^{\frak{F}}=G^{\frak{N}}N/N, $  and }

(ii) {\sl  $G^{\frak{N}}N=R^{\frak{N}}N$  provided $G=RN$}.

 {\bf Proof  of Proposition B.}  (i) Let $D=G^{{\mathfrak{N}}_{\sigma}}$. 
 First note that if  
$R$ is a normal subgroup of $G$, then 
$(G/R)^{{\mathfrak{N}}_{\sigma}}=DR/R$ by 
Lemmas 2.1 and   2.2 and  so from the $G$-isomorphism 
$DR/R\simeq D/(D\cap  R)$ we get that every chief factor of $G/R$ below
 $(G/R)^{{\mathfrak{N}}_{\sigma}}$ is simple if and only if every chief factor of $G$ 
between $D$ and $D\cap R$ is simple.
Therefore if  $G\in   {\mathfrak{U}}_{c\sigma}$, then $G/R\in   {\mathfrak{U}}_{c\sigma}$.   
 Hence the class  
${\mathfrak{U}}_{c\sigma}$  is closed under taking homomorphic images.

Now  we show that if  $G/R$, $G/N\in  {\mathfrak{U}}_{c\sigma}$,
 then $G/(R\cap N)  \in  {\mathfrak{U}}_{c\sigma}$.   We can assume without loss 
of generality that $R\cap N=1$. Since $G/R \in  {\mathfrak{U}}_{c\sigma}$,
every chief factor of $G$ 
between $D$ and $D\cap R$ is simple. Also, every chief factor of $G$ 
between $D$ and $D\cap N$ is simple. Now let $H/K$ be any chief factor of $G$ below $D\cap R$.  
Then $H\cap D\cap N=1$ and hence  from the $G$-isomorphism
 $$H(D\cap N)/K(D\cap N)\simeq H/(H\cap K(D\cap N))=H/K(H\cap D\cap N)=H/K$$ 
we get  that $H/K$ is simple since $D\cap N\leq K(D\cap N)\leq D$. On the 
other hand,    every chief factor of $G$ 
between $D$ and $D\cap R$ is also  simple. 
 Therefore the Jordan-H\"{o}lder 
theorem for groups with operators \cite[A, 3.2]{DH} implies that every chief factor of $G$ 
below $D$  is simple.  Hence   $G  \in  
{\mathfrak{U}}_{c\sigma}$, so  the class  
${\mathfrak{U}}_{c\sigma}$  is closed under taking subdirect products.

Finally, if  $H\trianglelefteq  G\in   {\mathfrak{U}}_{c\sigma}$, then from Lemmas 2.1 and  
2.2 and the isomorphism $$H/(H\cap D)\simeq HD/D \in  {\mathfrak{N}}_{\sigma}$$ we get that 
$H^{{\mathfrak{N}}_{\sigma}}\leq  H\cap D$ and 
so  every chief factor of $H$ 
below $H^{{\mathfrak{N}}_{\sigma}}$  is simple since every chief factor of $G$ 
below $D$ is simple. Hence
 $H\in   
{\mathfrak{U}}_{c\sigma}$,  so  the class  
${\mathfrak{U}}_{c\sigma}$  is closed under taking normal subgroups.

(ii)  See the proof of (i).

The proposition is proved.

{\bf Lemma 2.3.}   {\sl  Let $H/K$ be a non-abelian chief factor of $G$.
 If $H/K$ is simple, then $G/HC_{G}(H/K)$ is soluble.}

 {\bf Proof.}   Since   $C_{G}(H/K)/K=C_{G/K}(H/K)$, we can  assume without 
loss of generality that $K=1$. Then  $$G/C_{G}(H)\simeq V\leq 
\text{Aut}(H)$$ and $$H/(H\cap C_{G}(H) )\simeq HC_{G}(H)/C_{G}(H)\simeq \text{Inn}(H)$$
 since $C_{G}(H)\cap H=1.$ Hence 
$$G/HC_{G}(H)\simeq  (G/C_{G}(H))/(HC_{G}(H)/C_{G}(H))\simeq W\leq 
 \text{Aut}(H)/\text{Inn}(H).$$ From the validity of the Schreier 
conjecture, it follows that $G/HC_{G}(H/K)$ is soluble.           
The lemma is proved.

 {\bf Proof of Theorem D.}  First note that $D$ is characteristic in $N$ 
 and $R=D_{\mathfrak{S}}$ is a characteristic  subgroup of $D$,
  so both these subgroups are normal in $G$.

{\sl Necessity.}  In view of Proposition B(ii),  $G/G^{{\mathfrak{U}}_{\sigma}}$ is 
 $\sigma$-supersoluble  and 
$G^{{\mathfrak{U}}_{\sigma}}$ is contained in every normal subgroup $E$ of 
$G$ with $\sigma$-supersoluble quotient $G/E$. By Lemmas 2.1 and   2.2,
 $N/D=(G/N)^{{\mathfrak{N}}_{\sigma}}$.
 On the other hand, every 
chief factor of $G$ between $N$ and $D$ is abelian and so cyclic and hence  $G/D$ is  
$\sigma$-supersoluble.    Therefore $G^{{\mathfrak{U}}_{\sigma}}\leq D$.  
Moreover,  from Lemma   2.2 and Proposition B(ii) we also get that 
 $$N/G^{{\mathfrak{U}}_{\sigma}}=
(G/G^{{\mathfrak{U}}_{\sigma}})^{{\mathfrak{N}}_{\sigma}},$$ so every 
chief factor of $G$  between $N$ and  $G^{{\mathfrak{U}}_{\sigma}}$ is 
cyclic and hence $D\leq G^{{\mathfrak{U}}_{\sigma}}$. Thus
 $D=G^{{\mathfrak{U}}_{\sigma}}$, so if $D=1$, then $G$ is $\sigma$-supersoluble.

Now suppose that $D\ne 1$. We show that in this case $G$ has   a 
Robinson complex of the form   $(D,  Z(D); U_{1}, \ldots ,  U_{k})$, where $Z(D)=R$.  
It is clear that every chief factor of $G$ below $R$ is cyclic, so  
$G/C_{G}(R) $ is supersoluble by \cite[IV, 6.10]{DH}.  Hence $D=G^{{\mathfrak{U}}_{\sigma}}\leq C_{G}(R) $, 
so $R\leq   Z(D)\leq D_{\mathfrak{S}}=R$ and therefore we have  $ Z(D)=R$.

Now let $H/K$ be any chief factor of $G$ below $D$. Then $H\leq N$ and  so   in 
the case when $H/K$ is abelian, this factor  is cyclic, which implies that
 $D=G^{{\mathfrak{U}}_{\sigma}}\leq C_{G}(H/K)$. On the other hand, if $H/K$ is a non-abelian 
simple group, then Lemma 2.3 implies that $G/HC_{G}(H/K)$ is soluble. 
Then $$DHC_{G}(H/K)/HC_{G}(H/K)\simeq D/(D\cap HC_{G}(H/K))=D/HC_{D}(H/K)$$ 
is soluble, so  $D=HC_{D}(H/K)$ since $D$ is evidently perfect. Therefore, in both cases,  
every element of 
$D$ induces an inner automorphism on $H/K$. Therefore $D$ is 
quasinilpotent. Hence   in view of \cite[X, 13.6]{31},  $G$ has   a 
Robinson complex of the form   $(D,  Z(D), U_{1}, \ldots , U_{k})$.

{\sl Sufficiency.} From Conditions (i), (ii) and (iii), it follows that 
all  factors  below $N$ of any chief  series of $G$ passing through $N$ are simple. 
Therefore the Jordan-H\"{o}lder 
theorem for groups with operators \cite[A, 3.2]{DH} implies that every chief factor of $G$ 
below $N$  is simple.    Therefore $G$ is a  $\sigma$-$SC$-group. 

The theorem is proved. `

\section{Proofs of Proposition A and Theorem A}

Recall  that a subgroup $A$ of $G$ is called \emph{${\sigma}$-subnormal}
  in $G$ \cite{1} if   there is a subgroup chain  $$A=A_{0} \leq A_{1} \leq \cdots \leq
A_{n}=G$$  such that  either $A_{i-1} \trianglelefteq A_{i}$ or 
$A_{i}/(A_{i-1})_{A_{i}}$ is  ${\sigma}$-primary 
  for all $i=1, \ldots , n$.

{\bf Lemma 3.1} (See Remark 1.1 and [Proposition 2.6]{arivII}).  {\sl  $G$ is
  $\sigma$-nilpotent if and only if every 
subgroup of $G$ $\sigma$-subnormal in $G$.}

{\bf Lemma 3.2.} {\sl Let  $A$,  $K$ and $N$ be subgroups of  $G$.
 Suppose that   $A$
is $\sigma$-subnormal in $G$ and $N$ is normal in $G$.  }

(1) {\sl $A\cap K$    is  $\sigma$-subnormal in   $K$}.

(2) {\sl $AN/N$ is
$\sigma$-subnormal in $G/N$. }

(3) {\sl If $N\leq K$ and $K/N$ is
$\sigma$-subnormal in $G/N$, then $K$ is
$\sigma$-subnormal in $G.$}

(4) {\sl If $H\ne 1 $ is a Hall $\sigma _{i}$-subgroup of $G$  and $A$ is not  a
 $\sigma _{i}'$-group, then $A\cap H\ne 1$ is
 a Hall $\sigma _{i}$-subgroup of $A$. }

(5) {\sl If $A$ is a $\sigma _{i}$-group, then $A\leq O_{\sigma _{i}}(G)$. 
}

(6) {\sl  If $A$ is a Hall $\sigma _{i}$-subgroup of $G$, then $A$ is normal in $G$.}

(7) {\sl If  $|G:A|$ is a $\sigma _{i}$-number,  then  $O^{{\sigma _{i}}}(A)= O^{{\sigma _{i}}}(G)$.}

 (8)    {\sl  If $G$ is  $\sigma$-perfect, then $A$ is subnormal in $G$. }

 (9)    {\sl  $A^{{\frak{N}}_{\sigma}}$    is subnormal in $G$. }

{\bf Proof. }    Assume that this lemma   is false and let $G$ be a counterexample of
 minimal     order.   By hypothesis, there is a subgroup chain  $A=A_{0} \leq
A_{1} \leq \cdots \leq A_{r}=G$ such that
either $A_{i-1} \trianglelefteq A_{i}$ 
  or $A_{i}/(A_{i-1})_{A_{i}}$ is  $\sigma $-primary  for all $i=1, \ldots , r$. 
  Let   $M=A_{r-1}$.
  We can assume without loss of generality that $M\ne G$. 

(1)--(7) See Lemma 2.6 in \cite{1}.

 (8) $A$ is 
subnormal in $M$ by the choice of $G$. On the other hand, since $G$ is  $\sigma$-perfect,
 $G/M_{G}$  is not $\sigma$-primary. Hence $M$ is normal in $G$ and so $A$ 
is   subnormal in $G$.
 
(9)  $A$ is $\sigma$-subnormal in $AM_{G}\leq M$ by Part (1), so the 
choice of $G$ implies that $A^{{\frak{N}}_{\sigma}}$ is 
subnormal in $AM_{G}$.   Hence $G/M_{G}$ is a $\sigma 
_{i}$-group for some $i$, so $M_{G}A/M_{G}\simeq A/A\cap M_{G}$ is a $\sigma 
_{i}$-group.  Hence $A^{{\frak{N}}_{\sigma}}\leq  M_{G}$, so $A^{{\frak{N}}_{\sigma}}$ is subnormal
 in $M_{G}$ and hence  $A^{{\frak{N}}_{\sigma}}$ is subnormal
 in $G$.

Lemma is proved.

The following lemma, in fact, is a corollary of Theorem A and Lemmas 3.1 and 3.2(3). 

{\bf Lemma 3.3.}   {\sl  The following statements hold:}

(i) {\sl $G$ is a 
 $P\sigma T$-group if and only if every  $\sigma$-subnormal subgroup of 
$G$ is $\sigma$-permutable in $G$. }

(ii) {\sl 
If  $G$ is a 
 $P\sigma T$-group, then every   quotient $G/N$ of $G$ is also a   
$P\sigma T$-group. }

{\bf Lemma 3.4.}   {\sl 
Let  $A$ and $B$ be subgroups of $G$, where $A$ is
 $\sigma$-permutable  $G$. }

(1)  {\sl  If $A\leq B$ and $B$ is $\sigma$-subnormal in $G$, then  $A$ is
 $\sigma$-permutable  $B$}.

(2)  {\sl  Suppose that   $B$ is a $\sigma _{i}$-group. Then   $B$ is
 $\sigma$-permutable in  $G$  if and only if   $O^{\sigma _{i}}(G)  \leq 
N_{G}(B)$}.

 {\bf Proof. }  (1)  By hypothesis, $G$ possesses a complete Hall 
$\sigma$-set ${\cal H}=\{H_{1}, \ldots , H_{t}\}$. Then  ${\cal 
H}_{0}=\{H_{1}\cap B, \ldots , H_{t}\cap B\}$ is  a complete Hall 
$\sigma$-set of $B$ by Lemma 3.2(4). Moreover, for every  $x\in B$ and $H\in 
{\cal H}$  we have $AH^{x}=H^{x}A$,
 so $$AH^{x}\cap B=A(H^{x}\cap B)=A(H\cap B)^{x}=(H\cap B)^{x}A.$$   Hence  $A$ is
 $\sigma$-permutable in   $B$ by Remark 1.2. 

(2) See Lemma 3.1 in  \cite{1}.

The lemma is proved.

 {\bf Proof of Proposition A.}   Let ${\cal H}=\{H_{1}, \ldots , H_{t}\}$ and
  $N=G^{{\mathfrak{N}}_{\sigma}}$ be the $\sigma$-nilpotent residual of 
$G$. Then $D\leq N$.

(1) {\sl Statement (i) holds for $G$.}

Suppose that this  is false and let $G$ be a 
counterexample of minimal order.  If $D=1$, then $G$ is $\sigma$-soluble 
and so $G$ is a $\sigma$-$SC$-group  by Theorem B. Therefore $D\ne 1$.  
Let $R$ be  a minimal normal subgroup of $G$ contained in $D$. Then $G/R$ 
is a $P\sigma T$-group by Lemma 3.3(ii). Therefore the choice of $G$ 
implies that $G/R$ is a  $\sigma$-$SC$-group.    Since 
 $(G/R)^{{\mathfrak{N}}_{\sigma }}=N/R$  by Lemmas 2.1 and 2.2,
every 
chief  factor of $G/R$   below $N/R$ is simple. Hence   
 every 
chief  factor of $G$ between  $G^{{\mathfrak{N}}_{\sigma }}$   and $R$ is 
simple.   
Therefore, in view of the Jordan-H\"{o}lder 
theorem for groups with operators \cite[A, 3.2]{DH}, 
it is enough to show that 
 $R$ is   simple. Suppose 
that this is  false.    Let $L$ be 
a minimal normal subgroup of $R$.  Then $1  < L  < R$ and $L$ is 
$\sigma$-permutable in $G$ 
 by Lemma 3.3(i)   since 
$G$ is  a $P\sigma T$-group.   Moreover, $L_{G}=1$ and so $L$  is 
 $\sigma$-nilpotent by Theorem A.  Therefore  $R$ is a 
$\sigma _{i}$-group for some $i$, so for some $k$ we have $R\leq H_{k}$.
  Now let $V$ be a maximal subgroup of 
$R$. Then  $V$ is $\sigma$-subnormal in $G$, so $V$ is 
$\sigma$-permutable in $G$ and hence $$R\leq D\leq O^{\sigma _{i}}(G)\leq N_{G}(V)$$ by Lemma 3.4(2). 
Thus $R$ is nilpotent, so $R$ is a $p$-group for some $p\in \sigma _{i}$. 
Now let $V$ be a maximal subgroup of $R$ such that $V$ is normal in a 
Sylow $p$-subgroup of $P$ of $H_{k}$. By hypothesis,  $H_{k}$  is a $PST$-group and so $V$ 
is $S$-permutable in $H_{k}$ since it is subnormal in $H_{k}$. Then,
 by Lemma 3.4(2) (taking in the case $\sigma =\{\{2\}, \{3\}, \ldots 
\}$), we have  $H_{k}=PO^{p}(H_{k})\leq N_{G}(V)$.  
Therefore, in view of  Lemma 3.4(2),    
we have   $$G=H_{k}O^{\sigma _{i}}(G)\leq N_{G}(V).$$  Hence $V=1$ and so $|R|=p$, a 
contradiction. Thus we have (1).
 
(2)    {\sl Statement (ii) holds for $G$.}

  It is clear that $D$ is $\sigma$-perfect   and $G/D$ is  $\sigma$-soluble.
 In view of  Lemma 3.3(ii),  $G/D$    is a $P\sigma T$-group. 
 It is also cleat that  $D\leq G^{{\mathfrak{U}}_{\sigma}}$. On the other 
hand, $G/D$ is $\sigma$-supersoluble by Theorem B. Therefore     
$G^{{\mathfrak{U}}_{\sigma}}\leq D$ and so we  have $D=  
G^{{\mathfrak{U}}_{\sigma}}$. Hence we have (2).

(3)   {\sl Statement (iii) holds for $G$.}

Let $L$ be   a $\sigma$-soluble normal subgroup of $G$ and let $x$ be
 a $\sigma _{i}'$-element of $G$.  Let $V/L\leq O_{\sigma _{i}}(G/L)$.  
Then $V/L$ is  $\sigma$-subnormal in $G/L$, so $V/L$ is 
$\sigma$-permutable in $G/L$ by Lemma 3.3(i) since $G/L$ is a $P\sigma 
T$-group by  Lemma 3.3(ii). Therefore $$xL\in O^{\sigma _{i}}(G/L)\leq N_{G/L}(V/L)$$
 by Lemma 3.4(2).
Hence Statement (iii) holds for $G$.
 
The proposition is proved.

{\bf Lemma 3.5.}    {\sl Let $G$ be a non-$\sigma$-supersoluble $\sigma$-full 
  $\sigma$-$SC$-group and let    $(D, Z(D); U_{1}, \ldots , 
U_{k})$ be   a Robinson complex $G$, where $D=G^{{\mathfrak{U}}_{\sigma}}$. 
  Let $U$ be a
 non-$\sigma$-permutable $\sigma$-subnormal subgroup of $G$ of minimal 
order. Suppose that $S/Z(S)$ is $\sigma$-perfect.  Then:}

(i) {\sl If $US_{i}'/U_{i}'$ is $\sigma$-permutable  in  $G/U_{i}'$ for 
all $i$, then $U$ is $\sigma$-supersoluble.}

(ii) {\sl If $U$ is  $\sigma$-supersoluble and $UL/L$ is 
$\sigma$-permutable  in  $G/L$ for
  all non-trivial nilpotent  normal subgroups $L$ of $G$, then 
$U$ is a cyclic $p$-group for some prime $p$. }

{\bf Proof. }  Suppose that this lemma is false and let $G$ be a 
counterexample of minimal order.  By hypothesis, for some 
$i$ and for some  Hall $\sigma _{i}$-subgroup $H$ of $G$ we have $UH\ne 
HU$. 

(i) Assume that this is false.  Then $U\cap D\ne 1$ since $UD/D\simeq U/(U\cap D)$ is
 $\sigma$-supersoluble  by Proposition B(ii).  Moreover, Lemma 3.2(1)(2), 
 implies that $(U\cap D)Z(D)/Z(D)$ is 
$\sigma$-subnormal in $D/Z(D)$ and so  $(U\cap 
D)Z(D)/Z(D)$ is a non-trivial 
subnormal  subgroup of  $D/Z(D)$  by Lemma 3.2(8)
 since $D/Z(D)$ is $\sigma$-perfect by hypothesis. 
 Hence for some $i$ we 
have $$U_{i}/Z(D)\leq (U\cap 
D)Z(D)/Z(D),$$ so  $U_{i}\leq (U\cap 
D)Z(D).$ But then $$U_{i}'\leq  ((U\cap 
D)Z(D))'\leq U\cap D.$$  By hypothesis, $UU_{i}'/U_{i}'=U/U_{i}'$ is $\sigma$-permutable  in  
$G/U_{i}'$ and so
 $$UH/U_{i}'=(U/U_{i}')(HU_{i}'/U_{i}')=(HU_{i}'/U_{i}')(U/U_{i}')=HU/U_{i}'.$$ Hence $UH=HU$,
 a contradiction. Therefore Statement (i) holds.

(ii)  Let $N=
U^{{\mathfrak{N}}_{\sigma}}$.  Then $D$ is subnormal in $G$ by Lemma 
3.2(9). 
 Since $U$ is  $\sigma$-supersoluble by 
hypothesis, $N < U$. By Lemmas 2.1,  2.2 and 3.2(3),  every proper  subgroup $V$ of $U$ with 
$N\leq V$  is $\sigma$-subnormal in $G$, so the minimality of $U$ implies 
that $VH=HV$. Therefore, if $U$ has at least two distinct  maximal subgroups $V$ and $W$
 such that $N\leq V\cap W$, then $U=\langle V, W \rangle $ is permutes with $H$ by
 \cite[A, 1.6]{DH}, contrary to our assumption on $U$ and $H$. Hence $U/N$ 
is a cyclic $p$-group for some prime $p$.

 First assume that $p\in \sigma_{i} $. Lemma 3.2(4) implies that $H\cap U$  
is a Hall $\sigma _{i}$-subgroup  of $U$, so $U=N(H\cap U)=(H\cap U)N$. Hence  
$$UH=(H\cap U)NH=H(H\cap U)N=HU,$$ a contradiction. Thus  $p\in \sigma_{j}$ for some $j\ne i$.

Now we show that $U$ is  a $P\sigma T$-group. Let $V$ be a proper 
$\sigma$-subnormal subgroup of $U$.  Then $V$ is  
$\sigma$-subnormal in $G$  since $U$ is  
$\sigma$-subnormal in $G$.   The minimality of $U$ implies that $V$ is  
  $\sigma$-permutable   in $G$,  so $V$ is 
$\sigma$-permutable in $U$ by Lemma  3.4(1). Hence 
 $U$ is a $\sigma$-soluble $P\sigma T$-group by Lemma 3.3(i), so $N$ is 
abelian by Theorem B.

 Therefore $N$ is a $\sigma_{j}'$-group, so 
 $N\leq O=O_{\sigma_{j}'}(F(G))$ by Lemma 3.2(5) (taking in the case
 $\sigma =\{\{2\}, \{3\}, \ldots \}$).  By hypothesis, $OU/O$ permutes with 
$OH/O$. By Lemma 3.2(1)(2), $OU/O$ is $\sigma$-subnormal in 
$$(OU/O)(OH/O)=(OH/O)(OU/O)=OHU/O,$$
 where 
$OU/O\simeq U/U\cap O$
 is a $\sigma _{j}$-group and $OH/O\simeq H/H\cap O$ is a $\sigma _{i}$-group.
Hence $UO/O$ is normal in  $OHU/O$ by Lemma 3.2(6). Hence $H\leq N_{G}(OU)$
 $$H\leq N_{G}(O^{\sigma _{j}'}(OU))=N_{G}(O^{\sigma _{j}'}(U))$$ by Lemma 3.2(7)
 since
 $p \in \sigma _{j}$ implies that $O^{\sigma _{j}'}(U)=U$.
  But then $HU=UH$, a contradiction. 
Therefore Statement  (ii) holds.     

The lemma is proved.

{\bf Lemma 3.6.}   {\sl  Suppose that $G$ has a Robinson $\sigma$-complex $(D, Z(D); U_{1}, \ldots , 
U_{k})$, and let $N$ be a normal subgroup of $G$.  } 

(i) {\sl If $N=U_{i}'$, then  $$(D/N, Z(D/N)  =U_{i}/N; U_{1}N/N, \ldots , U_{i-1}N/N, 
U_{i+1}N/N, \ldots 
U_{k}N/N, U_{i}/N)$$ is  a Robinson $\sigma$-complex of $G/N$, where $U_{i}/U_{i}'\simeq
 Z(D)/(Z(D) \cap U_{i}')$. }

(ii) {\sl If $N$ is nilpotent, then  $$(DN/N, Z(DN/N)= Z(D)N/N; U_{1}N/N, \ldots , 
U_{k}N/N)$$ is  a Robinson $\sigma$-complex of $G/N$.}

 {\bf Proof. }  See Remark 1.6.8 in \cite{prod}.

{\bf Lemma 3.7} (See  Knyagina and  Monakhov \cite{knyag}). {\sl
Let $H$, $K$  and $N$ be pairwise permutable
subgroups of $G$ and  $H$ is a Hall subgroup of $G$. Then $$N\cap HK=(N\cap H)(N\cap K).$$}

{\bf Lemma 3.8.}   {\sl   If $G$  satisfies 
 $N_{\sigma _{i}}$, then $G/R$  satisfies 
 $N_{\sigma _{i}}$ for every    normal $\sigma$-soluble 
subgroup $R$ of $G$.}

 {\bf Proof. }   Let $N/R$ be a   normal $\sigma$-soluble
subgroup of $G/R$ and let $$(V/R)/(N/R)\leq O_{\sigma _{i}}((G/R)/(N/R)).$$ 
Then $N$ is a normal  $\sigma$-soluble
subgroup of $G$ and  $V/N\leq O_{\sigma _{i}}(G/N)$. Moreover,  
 for every $\sigma _{i}'$-element $xR\in G/R$ there is a   
$\sigma _{i}'$-element $y\in G$ such that $xR=yR$ and  so $yN\leq N_{G/N}(V/N)$, which implies that
 $$xR(N/R)\in N_{(G/R)/(N/R)}((V/R)/(N/R)).$$  
Hence  $G/R$  satisfies 
 $N_{\sigma _{i}}$, as required.

By the analogy with the notation   $\pi (n)$, we will write  $\sigma (n)$ to denote 
the set  $\{\sigma_{i} |\sigma_{i}\cap \pi (n)\ne 
 \emptyset  \}$;   $\sigma (G)=\sigma (|G|)$. 

{\bf Proof of Theorem C.}  First assume that $G$ is a $P\sigma T$-group and let
 $D=G^{{\mathfrak{S}}_{\sigma}}$  be the
 $\sigma$-soluble residual of $G$.  Then $D$ is clearly $\sigma$-perfect 
and, by 
Proposition A, $G$ is a $\sigma$-$SC$-group  and
 Statements (i) and (iii) hold for $G$. Moreover, Theorem B implies that 
$D$ coincides with the $\sigma$-supersoluble residual   $G^{{\mathfrak{U}}_{\sigma}}$ 
 of $G$ and if $D\ne 1$, then $G$ possesses a Robinson $\sigma$-complex of the form   $(D, 
 Z(D); U_{1},  \ldots U_{k})$.  Therefore the necessity of the condition 
of the theorem holds for $G$.

Now assume that $G$ has a normal $\sigma$-perfect subgroup $D$ and $D$ 
  satisfies    Conditions (i), (ii) and (iii). We show 
that   $G$ is a $P\sigma T$-group. Suppose that this is false 
 and let $G$ 
be a counterexample of minimal order.  
 Then $D\ne 1$ and $G$ has a $\sigma$-subnormal 
subgroup $U$ such that $UH\ne HU$ for some $i$ and some Hall $\sigma 
_{i}$-subgroup $H$ of $G$ and also every $\sigma$-subnormal 
subgroup $U_{0}$   of $G$   with $U_{0} < U$ is $\sigma$-permutable in 
$G$. Finally, note that $D=G^{{\mathfrak{U}}_{\sigma }}$ by Condition (i) and Theorem B.

(1) {\sl   $U$ is $\sigma$-supersoluble. 
  }

In view of Lemma 3.5(i), it is enough  to show that the hypothesis holds 
on    
$G/U_{i}'$ for all $i=1, \ldots , k$.  Let $N=U_{i}'$.  We can assume 
without loss of generality that $i=1$.   
 Then $$(D/N)^{{\mathfrak{N}}_{\sigma}}
 =D^{{\mathfrak{N}}_{\sigma}}N/N=D/N$$  by Lemmas 2.1 and  2.2, so $D/N$ is a normal
 $\sigma$-perfect subgroup
of $G/N$.  Moreover,   $(G/N)/(D/N)\simeq D/D$ is a $\sigma$-soluble $P\sigma T$-group.  
Now assume that $D/N\ne 1$.  Then, by                
 Lemma 3.6(i),   $$(D/N, Z(D/N); U_{2}N/N, \ldots U_{k}N/N)$$ is 
 a Robinson $\sigma$-complex of $G/N$, where $Z(D/N)=U_{1}/N$.    Moreover, if 
 $\{i_{1}, \ldots , i_{r}\}\subseteq \{2, \ldots , k\}$, where
 $2\leq r  < k$, then the quotients $G/N=G/U_{1}
'$ and $$(G/N) /(U_{i_{1}}N/N)'\cdots (U_{i_{r}}N/N)'=
(G/N)/(U_{i_{1}}'\cdots U_{i_{r}}'U_{1}'/N)\simeq G/U_{i_{1}}'\cdots U_{i_{r}}'U_{1}'$$
 satisfy
 $N_{\sigma _{l}}$ for all 
 $$\sigma _{l}\in  \sigma (U_{1}/N)= \sigma (Z(D/N)) \subseteq
 \sigma (Z(D)/(Z(D) \cap U_{1}')).$$  Therefore the hypothesis holds for 
$G/R$, so we have (1).

(2) {\sl $U$ is a cyclic $p$-group for some prime $p\in \sigma _{j}$, where $j\ne i$.}

First we show that $U$ is a cyclic $p$-group for some prime.
 In view of Claim (1) and Lemma 3.5(ii), it is enough  to show that the hypothesis holds 
on    
$G/N$ for every  normal nilpotent subgroup $N$ of $G$.   First note that 
$$(DN/N)^{{\mathfrak{N}}_{\sigma}}
 =D^{{\mathfrak{N}}_{\sigma}}N/N=DN/N$$  y Lemma 2.2(ii), so $D/N$
 is a normal $\sigma$-perfect subgroup  of $G/N$. Moreover, 
$$(DN/N,  Z(DN/N)=Z(D)N/N; U_{1}N/N, \dots , 
U_{k}N/N)$$   is  a Robinson $\sigma$-complex of $G/N$ by Lemma 3.6(ii).  Finally, if $V/N$ is a normal
 $\sigma$-soluble subgroup of $G/N$, then  $V$ is a normal
 $\sigma$-soluble subgroup of $G$ and so  for  
 $\{i_{1}, \ldots , i_{r}\}\subseteq \{1, \ldots , t\}$, where
 $1\leq r  < k$, the quotient   
  $G/N$ and, by Lemma 3.8, the quotient  $$(G/N) /(U_{i_{1}}N/N)'\cdots (U_{i_{r}}N/N)'=
(G/N) /(U_{i_{1}}'\cdots U_{i_{r}}'N/N)$$$$\simeq G/U_{i_{1}}'\cdots U_{i_{r}}'N\simeq
 (G/U_{i_{1}}'\cdots U_{i_{r}}')/(U_{i_{1}}'\cdots 
U_{i_{r}}'N/U_{i_{1}}'\cdots U_{i_{r}}') 
$$ satisfy
 $N_{\sigma _{l}}$ for all 
  for all $$\sigma _{l}\in \sigma (Z(DN/N))=\sigma (Z(D)N/N)\subseteq \sigma (Z(D))$$ since 
 $U_{i_{1}}'\cdots 
U_{i_{r}}'N/U_{i_{1}}'\cdots U_{i_{r}}'\simeq N/(N\cap U_{i_{1}}'\cdots 
U_{i_{r}}')$ is $\sigma $-soluble.

 Therefore   the hypothesis holds 
on $G/N$, so $U$ is a cyclic $p$-group for some prime $p\in \sigma _{j}$.  Finally, Lemma
 3.2(4) 
implies that in the case $i=j$ we have $U\leq H$, so $UH=H=HU$. Therefore 
 $j\ne i$.   Finally, again by Lemma 3.2(4), $U\leq O_{\sigma _{j}}(G)$.

 (3) $O_{\sigma _{j}}(G)\cap D=1$.

Suppose that $L=  O_{\sigma _{j}}(G)\cap D\ne 1$.   Then, since   $D/Z(D)$ is  $\sigma$-perfect, $L\leq Z(D)$  and  so $G$ satisfies 
 $N_{\sigma _{j}}$ by Condition (iii). Therefore $H\leq N_{G}(U)$ since 
$i\ne j$, $U\leq  O_{\sigma _{j}}(G)$ and $H$ is a ${\sigma _{i}}$-group. 
But then $HU=UH$, a contradiction. Hence we have (3).

{\sl Final contradiction for the sufficiency.}  By Lemma 3.2(2), $UD/D$ is 
$\sigma$-subnormal in $G/D$. On the other hand, $HD/D$ is a Hall $\sigma 
_{i}$-subgroup of $G/D$. Hence $$(UD/D)(HD/D)=(HD/D)(UD/D)=HUD/D$$ by Condition (i) and
 Lemma 3.3(i),
 so $HUD$ is a subgroup of $G$.
                  Therefore,   by Claims (2), (3) and Lemma 3.7, 
 $$UHD\cap HO_{\sigma _{j}}(G)
=UH(D\cap HO_{\sigma _{j}}(G)))=UH(D\cap H)(D\cap 
O_{\sigma _{j}}(G)))$$$$=UH(D\cap H)=UH$$ is a subgroup of $G$ and so $HU=UH$, a contradiction.

The theorem is proved.

\end{document}